\documentclass[12pt]{amsart}

\usepackage{latexsym}
\usepackage{amsmath}
\usepackage{amssymb}
\usepackage{amscd}
\usepackage{multicol}
\usepackage{mathrsfs}
\usepackage{xy}
\xyoption{all}

\newtheorem{theorem}{\bf Theorem}[section]

\newtheorem{proposition}[theorem]{\bf Proposition}
\newtheorem{lemma}[theorem]{\bf Lemma}

\newtheorem{definition-theorem}[theorem]{\bf Theorem-Definition}

\def\O{\mathbb{O}}
\def\H{\mathbb{H}}
\def\bR{\mathbb{R}}
\def\bC{\mathbb{C}}
\def\bZ{\mathbb{Z}}

\def\S{{\mathcal S}}
\def\L{{\mathcal L}}
\def\O{{\mathcal O}}
\def\o{{\mathbb O}}
\def\X{{\mathcal X}}
\def\Y{{\mathcal Y}}

  \def\w{\tau}
\def\quott({/\! /}

\def\h{\frak{h}}

\def\V{{\mathcal V}}
\def\X{{\mathcal X}}

\def\B{{\mathcal B}}

\title[Equivariant $K$-theory of quaternionic flag manifolds]
{Equivariant $K$-theory of quaternionic flag manifolds}
\author[A.-L. Mare ]{Augustin-Liviu Mare }
\author[M. Willems]{Matthieu Willems}
\date{\today}
\begin{document}
\address{{\rm Augustin-Liviu Mare},
Department of Mathematics and Statistics, University of
Regina, College West 307.14, Regina, Saskatchewan, S4S 0A2 Canada}
\email{mareal@math.uregina.ca}

\address{{\rm Matthieu Willems}, Department of Mathematics,
McGill University, 805 Sherbrooke Street West, Montr\'eal, Qu\'ebec, H3A 2K6 Canada}\email{matthieu.willems@polytechnique.org}

\begin{abstract}
We consider the manifold $Fl_n(\H)=Sp(n)/Sp(1)^n$ of
all complete flags in $\H^n$, where $\H$ is the skew-field of quaternions.
We study its equivariant complex $K$-theory rings with respect to
the action of two groups:  $Sp(1)^n$ and a certain canonical  subgroup $T=(S^1)^n$
(a maximal torus). For the first group action we obtain a Goresky-Kottwitz-MacPherson
type description. For the second one, we describe
the ring $K_T(Fl_n(\H))$ as a subring of $K_T(Sp(n)/T)$.
This ring is well known, since $Sp(n)/T$ is a complex flag variety.
\end{abstract}

\maketitle

\section{Introduction}
The quaternionic flag manifold $Fl_n(\H)$ is the space of all nested sequences
$$(V_{\nu})_{1\le \nu \le n}=V_1 \subset \ldots \subset V_n$$ where $V_{\nu}$ is a $\nu$-dimensional 
$\H$-vector subspace of $\H^n$, for all $1\le \nu \le n$, and
$\H$ is the skew-field of quaternions (by an $\H$-vector subspace we mean
a left $\H$-submodule). 
Let $Sp(n)$ denote the group of all $\H$-linear transformations of
$\H^n$ (that is, $n\times n$ matrices with coefficients in
$\H$) which preserve the canonical inner product on $\H^n$.
This group acts naturally on $Fl_n(\H)$.
In this paper we are particularly interested in the action on 
$Fl_n(\H)$ of the following two subgroups of $Sp(n)$: 
$$G=Sp(1)^n \  {\rm and } \ T=(S^1)^n. $$  
To put matters otherwise, $G$ is the group of all diagonal matrices in
$Sp(n)$ and $T\subset G$ consists of all such matrices with entries
in $\bC$, where $\bC$ is canonically embedded in $\H$
(as the set of all $a+bi$, where $a,b\in \bR$). It is worth mentioning that $G$ is actually 
the $Sp(n)$ stabilizer of the flag 
$(\H e_1\oplus \ldots\oplus  \H e_{\nu})_{1\le \nu \le n}$,
and since the action of $Sp(n)$ on $Fl_n(\H)$ is transitive, we can
identify
 $$Fl_n(\H) = Sp(n)/G.$$

We investigate the (complex, topological) equivariant $K$-theory rings
corresponding to the $T$ and $G$ actions. 
In general, the $G$-equivariant $K$-theory ring of any $G$-space $\X$ is
denoted by $K_G(\X)$ or $K_G^0(\X)$ (in this paper the first notation will be used in most
cases). By definition (see, for instance, \cite{Se}), it is the Grothendieck group of $G$-equivariant topological complex vector bundles over $\X$.  It is  a module over the ring 
$K_G({\rm pt.})=R[G]$, which is the representation ring of $G$.
We know that (see, for instance, \cite[Chapter 14, Section 6]{Hu})
 $$R[T]=R[(S^1)^n] = \bZ[x_1,x_1^{-1}, \ldots, x_n,x_n^{-1}]$$ 
and 
$$R[G]=R[Sp(1)^n]
=\bZ[X_1,\ldots,X_n].$$
Here $X_{\nu}=x_{\nu}+x_{\nu}^{-1},$ $ 1\le \nu \le n$,
are copies of  the (character of the) canonical representation of 
$Sp(1)=SU(2)$ on $\H=\bC^2$.

The flag manifold $Fl_n(\H)$ carries $n$ canonical (complex) vector bundles
$\V_1, \V_2, \ldots, \V_n$, of complex rank equal to 2, 4, $\ldots 2n$. 
The rank 2 quotient bundles $\L_{\nu}=\V_{\nu}/\V_{\nu-1},$ $1\le \nu \le n$ 
play an important role (by convention,  $\V_0$ is the rank 0 vector bundle).
Namely, we take into account that $T$ is a maximal torus in both $Sp(n)$ and $G$.
Moreover, the  Weyl group
 $W_G$ of $G$ is a normal subgroup of the Weyl group  
 $W_{Sp(n)}$ 
 and their quotient is
 $$W_{Sp(n)}/W_G \simeq \S_n,$$
 the symmetric group (see, for instance, \cite[Chapter 14, Section 4]{Hu}
 or Section \ref{roots}, below).
Results of \cite[Section 4]{Mc} (see also
Proposition \ref{themap1} of our paper), lead  to:
  \begin{align}\label{spin}
 K_G(Fl_n(\H))
 &\simeq K_{Sp(n)}(Fl_n(\H))\otimes_{R[Sp(n)]}R[G]\nonumber
 \\
 {}&\simeq R[G]\otimes_{R[Sp(n)]}R[G]
 \\
 {}&\simeq
\frac{ \bZ[[\L_1],\ldots,[\L_n], X_1,\ldots, X_n]}
{\langle \sigma_k([\L_1],\ldots,[\L_n]) -
\sigma_k(X_1,\ldots,X_n), 1\le k \le n\rangle}.\nonumber
\end{align}
Here $\sigma _k$ denotes the $k$-th symmetric polynomial in $n$ variables.
The ring $K_T(Fl_n(\H))$ 
is isomorphic to
$K_G(Fl_n(\H))\otimes_{R[G]}R[T]$ (see  \cite[Section 4]{Mc}
or Proposition \ref{thretwo}, below). Thus, we obtain
$$K_T(Fl_n(\H))
 \simeq
\frac{ \bZ[[\L_1],\ldots,[\L_n], x_1^{\pm1},\ldots, x_n^{\pm 1}]}
{\langle \sigma_k([\L_1],\ldots,[\L_n]) -
\sigma_k(x_1+x_1^{-1},\ldots,x_n+x_n^{-1}), 1\le k \le n\rangle}.
$$

In this paper we give alternative descriptions of the two rings 
 above.
The approaches will be different for $G$ and $T$, as follows.

The first main result  describes  $K_T(Fl_n(\H))$ as a subring of
the $T$-equivariant $K$-theory ring of the principal adjoint orbit $Sp(n)/T$.
The $T$-equivariant $K$-theory of principal adjoint orbits
(that is, complete complex flag varieties)
is  well understood, see, for instance,  \cite{Gr-Ra}, \cite{Ko-Ku},
 \cite{Le-Po}, \cite{Mc}, \cite{Li-Se}, \cite{Pi-Ra}, 
\cite{Wi1}, \cite{Wi2}.
If we identify $Sp(n)/T$ with the quotient of $Sp(2n,\bC)$ by  a Borel
subgroup, we deduce from the general theory (see, e.g.,
\cite[Lemma 4.9]{Ko-Ku}) that
$K_T(Sp(n)/T)$ has  a natural basis over $R[T]$, namely
the Schubert basis $\{[\O_w] \ : \ w\in W_{Sp(n)}\}.$
Like for any  flag variety, the  Weyl group $W_{Sp(n)}=N_{Sp(n)}(T)/T$ acts on $Sp(n)/T$ via
\begin{equation}\label{weyl}(nT).(gT)=gn^{-1}T,\end{equation}
for any $n\in N_{Sp(n)}(T)$ and $g\in Sp(n)$.
This action is $T$-equivariant. Therefore, by 
functoriality, it induces an action by ring homomorphisms on
$K_T(Sp(n)/T)$.
The canonical map $$\pi: Sp(n)/T\to Sp(n)/G=Fl_n(\H)$$
is $T$-equivariant too, hence it induces a homomorphism between the
$K_T$-rings, which we denote by $\pi^*_T$. We can now state the theorem.

\begin{theorem}\label{firstmain} The map
$\pi^*_T:K_T(Fl_n(\H)) \to K_T(Sp(n)/T)$ is injective.
Its image consists of all $W_G$-invariant elements
of $K_T(Sp(n)/T)$.
In this way, $K_T(Fl_n(\H))$ is the $R[T]$-subalgebra
of $K_T(Sp(n)/T)$ generated by all $[\O_w]$, where
$w\in W_{Sp(n)}$ is a maximal length representative of the
quotient $W_{Sp(n)}/W_G$.
\end{theorem}
 
 \noindent{\bf Remark.} A similar result holds for the general context of
 the ${\mathcal T}$-equivariant $K$-theories of ${\mathcal G}/{\mathcal P}$
 and ${\mathcal G}/{\mathcal B}$, where ${\mathcal G}$ is a complex semisimple
 Lie group, ${\mathcal P}$ a parabolic subgroup which contains
 the Borel subgroup  ${\mathcal B}$, and  ${\mathcal T}$ a maximal
 torus of ${\mathcal G}$ such that ${\mathcal T}\subset {\mathcal G}$
 (see, for instance, \cite[Corollary 3.20]{Ko-Ku}).
 However, Theorem \ref{firstmain} does not fit into this context, as
 $Fl_n(\H)$ is not a complex flag variety.


For $K_G(Fl_n(\H))$ we will prove the following   Goresky-Kottwitz-MacPherson
(shortly GKM) 
type description. Before stating it, we just mention that  the  $G$ fixed point set 
of $Fl_n(\H)$  can be identified with the symmetric group $\S_n$, as follows (see, for instance, \cite[Lemma 3.1]{Ma}):
\begin{equation}\label{fix}Fl_n(\H)^G
=\{(\H e_{\tau(1)}\oplus \ldots \oplus \H e_{\tau(\nu)})_{1\le \nu \le n}
\ : \ \tau\in \S_n\} = \S_n.\end{equation}
And here is the theorem.

\begin{theorem}\label{secondmain}
The ring homomorphism $K_G(Fl_n(\H)) \to \prod_{\tau\in \S_n}R[G]$ induced by the inclusion map $Fl_n(\H)^G\hookrightarrow  Fl_n(\H)$ is injective. Its image is
$$\{  (f_\tau) \in    \!\!\! \prod_{\tau \in \S_n} \!\!\! \ \mathbb{Z}[X_1, \ldots, X_n ] \ : \  f_\tau-f_{{(\mu, \nu)}\tau}  {\it \,\, is \,\, divisible \,\, by \,\,} X_\mu-X_\nu \,\,\, {\it 
for \ all } \  1 \leq \mu < \nu \leq n \}.
$$
\end{theorem}
Here ${(\mu, \nu)}$ denotes the transposition of $\mu$ and $\nu$, that is, the
 element of $\S_n$ which interchanges 
$\mu$ and $\nu$.

\noindent
{\bf Remarks.}  
1. A description of the integral cohomology ring of $Fl_n(\H)$ in terms of generators and relations has been obtained in
\cite{Bo} (see page 302). For the $G$-equivariant cohomology ring 
such a description has been obtained
in  \cite{Ma}. In both cases one obtains the same formulas as for the complex
flag manifold $Fl_n(\bC)$ (in the equivariant case the group acting on 
$Fl_n(\bC)$
is the standard maximal torus $T$ of the unitary group $U(n)$).
In the present paper we show that the same similarity can be noted for $K$-theory. 
Not only has  the ring $K_G(Fl_n(\H))$  the same presentation as $K_T(Fl_n(\bC))$
(see equation (\ref{spin})), but  also the same GKM description holds true (see Theorem
\ref{secondmain} and compare with \cite[Theorem 1.6]{Mc} for $G=U(n)$).
Another space for  which we have the same  analogy  with the complex flag manifold 
at the level of equivariant cohomology and $K$-theory is the octonionic flag manifold 
$Fl(\o)$  (see the recent paper \cite{Ma-Wi}). It would be interesting to find more
examples of spaces with group actions for which the equivariant $K$-theory has the
same features as $Fl_n(\H)$ and $Fl(\o)$. 
What makes these two spaces special is as follows:
First,
they are  homogeneous, of the form ${\mathcal G}/{\mathcal H}$ where 
${\mathcal G}$ is a compact
Lie group and ${\mathcal H}$ a closed subgroup of the same rank as ${\mathcal G}$; 
the group action is the one of
${\mathcal H}$, by  multiplication from the left. Second, a maximal torus 
${\mathcal T}$ of ${\mathcal H}$ has  the same fixed point set as ${\mathcal H}$ itself.
Third, ${\mathcal G}/{\mathcal H}$ admits a cell decomposition such that
each cell is ${\mathcal T}$-invariant, homeomorphic to $\bC^m$ for some $m\ge 0$,
and the action of ${\mathcal T}$ on it is complex linear. 

2. It is  worth noticing that all the previously known
GKM type descriptions of equivariant $K$-theory are for
actions of tori (for instance, any subvariety of a complex projective space
which is preserved
by a linear torus action, see, for instance, \cite[Appendix A]{Ro}).  
Theorem \ref{secondmain} is a non-abelian version of this general result.

3. A GKM  description exists for $K_T(Fl_n(\H))$ as well
(see Proposition \ref{gkmc}, below).
 It can also be deduced from the fact that $Fl_n(\H)$ has a cell decomposition
 whose cells are complex vector spaces, the torus $T$ leaving them invariant and
 acting on them complex linearly in a very explicit way (see  Section \ref{second}, below). Thus, one can 
apply the main result of 
\cite{Ha-He-Ho}: the main ingredient is the calculation of the Euler class
in $K_T$ for any  cell and the observation that this is not a zero-divisor
in $K_T({\rm pt.})=R[T]$. We will not present the details.
Our proof of Proposition \ref{gkmc} goes along  different lines, using the GKM description of
$K_T(Sp(n)/T)$.
Concerning 
 $K_G(Fl_n(\H))$,  we do not know if the GKM description  given in Theorem \ref{secondmain} is a  direct consequence of  \cite{Ha-He-Ho}. Even though the cells mentioned above are $G$-invariant
and the action of $G$ on cells is also very explicit (see again Section \ref{second}), 
this action is $\bR$-linear without being $\bC$-linear. It seems difficult to find  a  way to
compute the corresponding Euler classes in the ring $K_G({\rm pt.})=R[G]$ of 
{\it complex} representations of $G$.  

\noindent{\bf Acknowledgements.} We would like to thank the referees  for several valuable suggestions.

\section{The roots of $Sp(n)$}\label{roots}

In this section we collect some background material concerning
the roots of $Sp(n)$ and other related objects.
The details can be found  for instance in 
\cite[Section 16.1]{Fu-Ha}.

Let $\mathfrak{sp}(2n, \bC)$ be the  complexified Lie algebra of 
$Sp(n)$.
It consists of all complex square  matrices of the form $$\left( \begin{array}{cc}
 a & b \\
c & -a^t 
\end{array} \right),$$ where $a,b,c$ are $n\times n$ complex matrices
with $b^t = b$, $c^t = c$.
The elements of the complexified Lie algebra of $T$, call it
$\h$,  
are block matrices  as above with $a$ diagonal and $b=c=0$. 
A linear basis of $\h$ over $\bC$ consists of the matrices 
 $e_{\nu,\nu}-e_{\nu+n,\nu+n},$
 $1 \leq \nu \leq n$, where $\{ e_{\mu,\nu}\}_{1 \leq \mu,\nu \leq 2n}$ is the canonical basis of the space of complex $2n\times 2n$ matrices. 
 We denote by $\{ L^{\nu} \ : \ 1\le \nu \le n\}$ the  corresponding dual basis of $\mathfrak{h}^*$. We have as follows:
 \begin{itemize}
\item The set of roots is $$\Delta = \{ \pm L^{\mu}\pm L^{\nu}  \ : \  1 \leq \mu < \nu \leq n \} \cup \{ \pm 2L^{\nu}, 1 \leq \nu \leq n \}.$$
\item The set of positive roots is $$\Delta^+ = \{ L^{\mu}\pm L^{\nu} , 1 \leq \mu < \nu \leq n \} \cup \{ 2L^{\nu}, 1 \leq \nu \leq n \}.$$
\item A simple root system is 
$$\Pi = \{ \alpha_1 = L^1 - L^2 ,   \alpha_2 = L^2 - L^3, \ldots,  \alpha_{n-1} = L^{n-1} - L^{n}, \alpha_n = 2L^n\}.$$
\item The weight lattice $\mathfrak{h}_{\mathbb{Z}}^*$ is the $\mathbb{Z}$-module generated by $L^1, L^2, \ldots, L^n$.
\item The weight lattice $\mathfrak{h}_{\mathbb{Z}}^*$ is canonically isomorphic to $R[T] $, and $L^{\nu}$ corresponds to $x_{\nu}$ in this isomorphism,
$1\le \nu \le n$. More generally, for any $ \lambda \in \mathfrak{h}_{\mathbb{Z}}^*$, we will denote  the corresponding character in $R[T]$  by $e^{\lambda}$.
\item The Weyl group $W_{Sp(n)}=N_{Sp(n)}(T)/T $ is generated by the simple reflections $s_{\nu}$ corresponding to $\alpha_{\nu}$, for $1 \leq \nu \leq n$
(as usual,  we denote by $s_{\alpha}$ the reflection corresponding to the
positive root $\alpha$).
Concretely, $W_{Sp(n)}$ consists of all linear  automorphisms $\eta$ of 
$\mathfrak{h}^*$ such that for any $1\leq \nu \leq n$, there exists $1\leq \mu \leq n$ 
such that $\eta(L^{\nu})=\pm L^{\mu}$. This means that $W_{Sp(n)}$  is the semi-direct product of the symmetric group $\S_n$ of permutations of the set $\{ L^{\nu} \ : \  1 \leq \nu \leq n \}$ and the group $\{ -1,1 \}^n$ of sign changes (here
$\S_n$ acts on $\{-1,1\}^n$ by permuting the entries of an $n$-tuple). 
For $1 \leq \nu \leq n-1$, the reflection $s_{\nu}$ is the transposition $(\nu,\nu+1)$. The reflection $s_n$ sends $L^n$ to $-L^n$, and $L^{\nu}$ to itself, for $1 \leq \nu < n$.
More generally, for $1 \leq \mu < \nu \leq n$, the reflection corresponding to 
the root $L^{\mu}-L^{\nu}$ is the transposition $(\mu,\nu)$, whereas the reflection corresponding to $L^{\mu}+L^{\nu}$ sends $L^{\mu}$ to $-L^{\nu}$, $L^{\nu}$ to 
$-L^{\mu}$, and leaves $L^{\kappa}$ unchanged, for $\kappa\notin
\{\mu,\nu\}$. The reflection corresponding to the root $2L^{\nu}$ sends 
$L^{\nu}$ to $-L^{\nu}$ and leaves $L^{\mu}$ unchanged  for $\mu \neq \nu$.
\item The Weyl group $W_{G}: = N_G(T)/T$ is the subgroup of
$W_{Sp(n)}$  generated by $s_{2L^{\nu}}$,  $1 \leq \nu \leq n$. 
It is isomorphic to $\{ -1,1 \}^n$.
\item The quotient $ W_{Sp(n)} / W_{G}$ is isomorphic to the group of permutations of the set $\{ L^\nu   \ : \ 1 \leq \nu \leq n \}$, which  is the symmetric group $\S_n$.
For any pair $\mu, \nu$ such that $1 \le \mu <\nu \le n$, the cosets
$s_{L^\mu -L^\nu}W_G$ and $s_{L^\mu+L^\nu}W_G$ are equal and 
are mapped by the  isomorphism above to
the transposition $(\mu,\nu)$.
\end{itemize}

\section{The cell decomposition}\label{second}
 
 In this section we describe the Schubert cell decomposition of $Fl_n(\H)$.
 We will be especially interested in the action of $G$ on the cells.

 The group $GL_n(\H)$ of all invertible $n\times n$ matrices with entries in 
 $\H$ acts linearly on $\H^n$.
 More precisely,  $\H^n$ is regarded as a left $\H$-module 
 and the action of $GL_n(\H)$ is given by:
 $$gh=h\cdot g^*$$ 
 for any $g\in GL_n(\H)$ and any $h\in \H^n$. Here
 $\cdot$ denotes the matrix multiplication and $g^*$ the transposed conjugate of $g$. 
The group $GL_n(\H)$ acts  on $Fl_n(\H)$ by 
 $$g(V_{\nu})_{1\le \nu \le n}=(gV_{\nu})_{1\le \nu \le n},$$
 for any $g\in GL_n(\H)$ and any $(V_{\nu})_{1\le \nu \le n}\in Fl_n(\H)$.
 This group action is transitive and the stabilizer of the
 flag $(\H e_1\oplus \ldots \oplus \H e_{\nu})_{1\le \nu \le n}$ is the 
 group $B$ consisting of all upper triangular matrices with entries in
 $\H$.   
 In this way we obtain the identification
 $$Fl_n(\H)=GL_n(\H)/B.$$
 The following result is a direct consequence of the Bruhat decomposition of
 $GL_n(\H)$ (see 
 \cite[Section 19, Theorem 1]{Dr}):
 
 \begin{proposition}\label{dec} Any $g\in GL_n(\H)$ can be  written as
 $g=up_\w b$, where:
 
 \begin{itemize}
 \item[1.]  
   $\w\in {\mathcal S}_n$ and 
 $p_\w$ denotes the matrix $(\delta_{\mu, \w(\nu)})_{1\le \nu,\mu\le n}$,
 where $\delta$ is the Kroenecker delta. 
 \item[2.]  $b\in B$
 \item[3.]   both $u$  and $(p_\w up_\w^{-1})^t$ are upper triangular with all entries on the diagonal  equal to 1 (the superscript $t$ indicates the matrix transposed).
 \end{itemize}
 
 Moreover, the matrices $p_\w$ and $u$ with  properties 1 and 3 above 
 are uniquely determined by $g$.
 \end{proposition} 
 
 We deduce that
\begin{equation}\label{cell}GL_n(\H)/B=\bigsqcup_{\w\in {\mathcal S}_n} {\mathfrak U}_\w p_\w B/B,\end{equation}
 where ${\mathfrak U}_\w$ denotes the set of  all $n\times n$ matrices $u$ with
 entries in $\H$   such that both  $u$ and $(p_\w up_\w^{-1})^t$ are upper triangular with all entries on the diagonal  equal to 1. The canonical map ${\mathfrak U}_\w
 \to {\mathfrak U}_\w p_\w B/B$ is a homeomorphism. Indeed, this map is
 continuous and bijective, by Proposition \ref{dec}.  Its inverse is continuous too,
 because 
the map ${\mathfrak U}_\w p_\w B\to {\mathfrak U}_\w$ which assigns 
to $g=up_\w b$ the first factor $u$  is continuous, as we can see from its explicit description  in the proof of   \cite[Section 19, Theorem 1]{Dr}. 
Now it is an easy exercise to see that an $n\times n$ matrix
$u=(u_{\mu\nu})_{1\le \mu,\nu\le n}$ is in   ${\mathfrak U}_\w$ 
if and only if the diagonal entries are equal to 1 and the
others are equal to 0, except for those $u_{\mu\nu}$ with
$\mu<\nu$ and $\w (\mu)>\w (\nu)$ (the key point is the formula
$p_\w up_\w^{-1}=(u_{\w(\mu)\w(\nu)})_{1\le \mu,\nu\le n}$). This implies that ${\mathfrak U}_\w$ can be
identified with $\H^{\ell(\w)}$, where $\ell(\w)$ denotes the number of inversions 
of the permutation $\w$. Consequently, for any $\w\in \S_n$ the  element 
$$C_\w={\mathfrak U}_\w p_\w B/B$$ of the
decomposition (\ref{cell}) is homeomorphic to a  cell of (real) dimension $4\ell(\w)$. 
We call it a {\it Bruhat cell}. 

The group $G$ acts on $GL_n(\H)/B$ by left multiplication. We claim that this
action leaves any Bruhat cell $C_\w={\mathfrak U}_\w p_\w B/B$ invariant. 
To justify this, take $\gamma={\rm Diag}(\gamma_1,\ldots,\gamma_n)$  in  $G$.
We have
$$\gamma p_\w=
{\rm Diag}(\gamma_1,\ldots,\gamma_n)p_\w=  p_\w{\rm Diag}(\gamma_{\w^{-1}(1)},\ldots,\gamma_{\w^{-1}(n)}).$$
This implies that if $u\in {\mathfrak U}_\w$, then
$$\gamma u p_\w B=\gamma u  \gamma^{-1} \gamma p_\w B
=\gamma u \gamma^{-1} p_\w B.$$
We  notice that if  $u=(u_{\mu\nu})_{1\le \mu,\nu\le n}$, then
$\gamma u \gamma^{-1}
=(\gamma_{\mu}u_{\mu\nu}\gamma_{\nu}^{-1})_{1\le \mu,\nu\le n}$. Thus, if $u$ is in $ {\mathfrak U}_\w$, then $\gamma u\gamma^{-1}$ is in 
${\mathfrak U}_\w$ as well.

We summarize our previous discussion as follows:

 \begin{proposition}\label{cells}  
 The Bruhat cell decomposition of the quaternionic flag manifold $Fl_n(\H)$ is 
 $ Fl_n(\H) =\bigsqcup_{\w\in \S_n} C_\w.$
The cell $C_\w$  has real dimension $4\ell(\w)$ and is $G$-invariant, being in fact  $G$-equivariantly homeomorphic to
$\bigoplus_{(\mu,\nu)} \H_{\mu\nu}$.
Here the  sum runs over all pairs $(\mu,\nu)$ with $1\le \mu <\nu\le n$ such that
$\w(\mu)>\w(\nu)$, and $\H_{\mu\nu}$ is a copy of $\H$.
The action of  $G$ on $\H_{\mu\nu}$ is
\begin{equation}\label{ga}(\gamma_1,\ldots,\gamma_n).h
=\gamma_\mu h\gamma_\nu^{-1},\end{equation}
for all $(\gamma_1,\ldots,\gamma_n)\in G$ and $h\in \H$.
 \end{proposition}

 \noindent{\bf Remark.} Let us identify 
 $$\H=\bC \oplus j\bC=\bC^2.$$
 It is an easy exercise to see that if in equation (\ref{ga}) we take
 $(\gamma_1,\ldots, \gamma_n)\in T$,
 the resulting transformation of $\H$ is $\bC$-linear.
 In other words, $T$ acts complex linearly on each cell $C_\w$.
 However, if $ (\gamma_1,\ldots,\gamma_n)$ is in $G$ but not
 in $T$, the transformation is in general not $\bC$-linear.

There are two alternative presentations of the cell $C_\tau$, which are
given in what follows. 

\begin{lemma}\label{albe} We have
$$C_\tau = Bp_\tau B/B.$$
\end{lemma}

\begin{proof} We have
$$GL_n(\H)=\bigsqcup_{\w\in {\mathcal S}_n} {\mathfrak U}_\w p_\w B$$
and ${\mathfrak U}_\w \subset B$ for all $\tau \in \S_n$.
Thus, it is sufficient to prove that if $\tau_1,\tau_2\in \S_n$,
$\tau_1 \neq \tau_2$, then 
$(Bp_{\tau_1} B)\cap (Bp_{\tau_2}B) =\emptyset$.
This can be proved by using the same arguments as in the proof of
\cite[Chapter III, Proposition 4.6]{Hi}.
\end{proof}

We can also describe $C_\tau$ by using the actual definition of
$Fl_n(\H)$, as the set of all flags in $\H^n$. 
Our model is the presentation of the Bruhat 
cells in $Fl_n(\bC)$, as given, for instance, in  \cite[Chapter III, Section 4]{Hi}.
First,  to each $r$-dimensional linear subspace $V$ of $\H^n$ we assign
the set $s(V)=\{m_1,\ldots,m_r\}$ where $m_1<\ldots <m_r$ are determined by
$V\cap \H^{m_t-1}\neq V\cap \H^{m_t}$, $1\le t \le r$
(here $\H^m$ denotes  $\H e_1 \oplus \ldots \oplus \H e_m$, for all
$1\le m \le n$, and $\H^0=\{0\}$).   Note that
if $V$ and $V'$ are subspaces such that $V\subset V'$, then
$s(V)\subset s(V')$.
To the
flag $V_\bullet=(V_\nu)_{1\le \nu\le n}$ we assign
the permutation $\tau=\tau^{V_\bullet}$ which is defined  recursively as follows:
\begin{itemize}
\item  $\{\tau(1)\}=s(V_1)$.
\item if $\tau(1),\ldots,\tau(k)$ are known, we set 
$\{\tau(k+1)\}=s(V_{k+1})\setminus s(V_k)$ 
\end{itemize}

\begin{lemma} We have $$C_\tau=\{V_\bullet \in Fl_n(\H)\ : \ 
\tau^{V_\bullet}=\tau\}.$$ 
\end{lemma}

\begin{proof} We are using the identification $GL_n(\H)/B=Fl_n(\H)$
given by $$gB = g \H_\bullet,$$ for any $g\in GL_n(\H)$.
Here $\H_\bullet$ denotes the flag 
$(\H e_1 \oplus \ldots \oplus \H e_\nu)_{1\le \nu \le n}$. 
In this way,  $p_\tau B$ is the same as the flag
$(\H e_{\tau(1)} \oplus \ldots \oplus \H e_{\tau(\nu)})_{1\le \nu \le n}$,
which we denote by $\tau \H_\bullet$.
The $B$-orbit of this flag is just $C_\tau$ (see Lemma \ref{albe}).
The lemma is a straightforward consequence of the following two facts:
\begin{align*}
{}& \tau^{\tau \H_\bullet}=\tau\\
{}& \tau^{b V_\bullet}=\tau^{V_\bullet} {\rm \ for \ all \ } b\in B 
 \ {\rm and \ all \ }V_\bullet \in Fl_n(\H).
\end{align*} 
\end{proof} 
In the same way as in \cite[p. 122]{Hi}, we see that the closure of $C_\w$
in $Fl_n(\H)$ consists
of all flags $V_\bullet$ such that $\tau^{V_\bullet}\preceq \tau$. Here
$\preceq$ denotes the Bruhat ordering on the symmetric group $\S_n$
(see \cite[Chapter I, Section 6]{Hi}). Note that if $\tau_1\preceq \tau_2$ then
$\ell (\tau_1)\le \ell(\tau_2)$.
\begin{proposition}\label{closure} The closure of $C_\w$ in $Fl_n(\H)$ can be expressed as follows:
$$\overline{C_\tau}=\bigsqcup_{\tau'\preceq \tau} C_{\tau'}.$$
Any of the cells $C_{\tau'}$ above,  with $\tau'\neq \tau$, has dimension strictly less than the dimension of $C_\tau$.
\end{proposition}

\section{$T$-equivariant $K$-theory}\label{prev}

Throughout this section we will use the notations
$$\X=Fl_n(\H)=Sp(n)/G\ {\rm and } \ \Y=Sp(n)/T.$$

Our main goal here is to prove Theorem \ref{firstmain}.
We will use the injectivity of the restriction to fixed points, which
is the content of the next proposition. We  first note  that
 the $T$ and
$G$ fixed points of $Fl_n(\H)$ are the same    (see \cite[Lemma 3.1]{Ma}).
By equation (\ref{fix}) we have
$$\X^G =\X^T=\S_n.$$
Let $i : \S_n \rightarrow \X$ be the inclusion map and $i^*_T : K_T(\X) \rightarrow 
\prod_{\w \in \S_n} R[T]$ the corresponding ring homomorphism.

\begin{proposition}\label{inject}

(a) The $T$-equivariant $K$-theory of $\X$ is a free $R[T]$-module of rank $n!$. 

(b) The restriction to fixed points $i^*_T : K_T(\X) \rightarrow \prod_{\w \in \S_n} R[T]$ is injective.

\end{proposition}

\begin{proof}

(a) Let $\X=\bigsqcup_{\w \in \mathcal{S}_n} C_\w$ be the cell decomposition of $\X$ (see Proposition \ref{cells}). Each cell $C_\w$ is $T$-equivariantly homeomorphic to $\mathbb{C}^{2\ell(\w)}$. Notice that
$$\ell (\w)\le \frac{n(n-1)}{2},$$
for all $\w\in \S_n$.  For any integer number
$\nu$ with $0 \leq \nu \leq \frac{n(n-1)}{2}$, we set 
$$\X_{\nu}=\bigsqcup_{\w \in \mathcal{S}_n, \ell (\w) \leq \nu} C_\w.$$ 
By Proposition \ref{closure}, this is a closed subspace of $\X$. 
We  prove by induction that for all $0 \leq \nu \leq \frac{n(n-1)}{2}$
we have:
\begin{itemize}
\item $K_T(\X_{\nu})=K_T^0(\X_{\nu})$ is a free $R[T]$-module of rank equal to the number of cells in $\X_{\nu}$
\item   $K_T^{-1}(\X_\nu)=\{0\}$. 
\end{itemize}
The definition of the functors $K_T^{-1}$ and $K_T^1$ 
 can be found for instance in  \cite{Se}. We have  
\begin{equation}\label{pont}K_T^{-1}({\rm pt.})=K_T^1({\rm pt.})=\{0\}.\end{equation}
This implies our claim for $\nu=0$, since  $\X_0$ consists of only one point
(of course we have $K_T({\rm pt.})=R[T]$). 

Let us assume that the claim is true for $\nu$. We will  prove it for $\nu+1$. 
We consider the space 
\begin{equation}\label{xnu}  \X_{\nu+1}\setminus \X_\nu =\bigsqcup_{\w\in \S_n, \ell(\w)=\nu+1}C_\w.\end{equation}
We consider now the exact sequence of the pair $(\X_{\nu+1},\X_\nu)$.
Since $\X_\nu$ is closed in $\X_{\nu+1}$, we deduce that
$K^*_T(\X_{\nu+1},\X_\nu)=K^*_T(\X_{\nu+1}\setminus\X_\nu)$ and obtain the
following sequence
(see \cite[Section 2]{Se}):
\begin{align*}{}& K_T^{-1}(\X_{\nu+1} \setminus \X_\nu) \rightarrow  K_T^{-1}(\X_{\nu+1}) \rightarrow K_T^{-1}(\X_{\nu}) \rightarrow 
\\
 {}& \rightarrow K_T^{0}(\X_{\nu+1} \setminus \X_{\nu}) \rightarrow  
 K_T^{0}(\X_{\nu+1}) \rightarrow K_T^{0}(\X_{\nu}) \rightarrow 
K_T^{1}(\X_{\nu+1} \setminus \X_{\nu}).
\end{align*}
By the induction hypothesis, $K_T^{-1}(\X_{\nu})=0$, and $K_T^0(\X_{\nu})$ is a free $R[T]$-module of rank equal to the number of cells in $\X_{\nu}$. 
Each of the cells $C_\w$ in the  union given by equation (\ref{xnu}) is a connected component
of $\X_{\nu+1}\setminus \X_\nu$:
indeed, by Proposition \ref{closure}, any such $C_\w$ is a closed subspace of
$  \X_{\nu+1}\setminus \X_\nu$, hence it is open as well (as its complement is closed).
We obtain
 $$K_T^{1}(\X_{\nu+1} \setminus \X_{\nu}) = K_T^{-1}(\X_{\nu+1} \setminus \X_{\nu})=\{0\},$$ 
 where we have used the Thom isomorphism for each of the cells $C_\w$ in the
 union given by (\ref{xnu}) (recall that $C_\w$ is  $T$-equivariantly homeomorphic  to 
 $\bC^{2\ell(\w)}$, the action of $T$ being complex linear, as mentioned in 
 the remark following Proposition \ref{cells}).
 We also have 
  $$K_T^0(\X_{\nu +1} \setminus \X_{\nu})  \simeq \bigoplus_{\w \in \S_n, \ell(\w) = \nu +1} K_T^0(C_\w)=
\bigoplus_{\w \in \S_n, \ell(\w) = \nu +1}R[T] .$$
In other words, $K_T^0(\X_{\nu +1} \setminus \X_{\nu}) $ 
is a free $R[T]$-module of rank equal to the number of cells  in $\X_{\nu+1} \setminus \X_{\nu}$. The desired conclusion follows.
 
 Finally, for $\nu=\frac{n(n-1)}{2}$, we obtain point (a) of the proposition.

Point (b) is a consequence of point (a). Indeed, let $Q[T]$ denote
the fraction field of $R[T]$. According to the localization theorem
(see \cite{Se}), the homomorphism $K_T(\X)\otimes_{R[T]}Q[T]
\rightarrow \prod_{\w \in \mathcal{S}_n}Q[T]$ induced by $i_T^*$ is
an isomorphism. Moreover, since $K_T(\X)$ is a free $R[T]$-module,
the canonical map  $K_{T}(\X) \rightarrow K_{T}(\X)\otimes_{R[T]}Q[T]
$ is an embedding. Since the  following  diagram
is commutative, we deduce  that   $i_{T}^*$ is injective.
$$\xymatrix{
    K_{T}(\X)\ar[dd]^{i_T^*} \ar@{^{(}->}[rr] & &
   K_{T}(\X)\otimes_{R[T]}Q[T]\ar[dd]^{\simeq} \\
    \\
   \prod_{ \mathcal{S}_n} R[T] \ar[rr] & &
  \prod_{ \mathcal{S}_n} Q[T] }
 \\
$$
\end{proof}

\noindent{\bf Remark.} The arguments used in the proof 
are standard: see \cite[proof of Lemma 2.2]{Mc} or
\cite[proof of Theorem 3.13]{Ko-Ku}.

Let $\rho : G\times \bC^m \to \bC^m$ be a complex representation of
$G$. The quotient space
\begin{equation}\label{vro}V_{\rho}=Sp(n)\times \bC^m/((k,v)\sim (kg,\rho(g^{-1})v)
{\rm \ for \ all \ } k\in Sp(n), g\in G, v\in \bC^m)\end{equation}
has a natural  structure of a vector bundle over $Sp(n)/G=Fl_n(\H)$.
Moreover,  there is a natural action of $G$ (consequently, also of $T$) on $V_{\rho}$ given by 
\begin{equation}\label{gkv}g(k,v)=(gk,v), \, {\rm \ for \ all \ }  g \in G, k \in Sp(n), { \rm and \, } v \in \mathbb{C}^m.
\end{equation}
The $R[T]$-linear extension of the assignment $\rho \mapsto V_{\rho}$
gives the  homomorphism $ \kappa_T : R[T] \otimes R[G] \rightarrow K_T(\X).$ It is easy to check that if $\rho$ is a representation of
$Sp(n)$, then $V_{\rho}$ is isomorphic to the trivial bundle
$(Sp(n)/G) \times \bC^m$. In other words, if  $\chi \in R[Sp(n)]$, then we have $\kappa_T(\chi\otimes 1-1 \otimes \chi)=0$. 
Consequently, we obtain a homomorphism
$$\overline{\kappa}_T : R[T] \otimes_{R[Sp(n)]} R[G]  \rightarrow K_T(\X).
$$

We are now ready to prove the first part of Theorem \ref{firstmain},
concerning the image of the map
$$\pi^*_T:K_T(\X)\to K_T(\Y).$$

\begin{proposition}\label{thretwo}

The homomorphism $\pi_T^*$ in Theorem \ref{firstmain} is injective and its image is equal to $K_T(\Y)^{W_{G}}$. Moreover, $\overline{\kappa}_T$ is a ring isomorphism.
\end{proposition}

\begin{proof}
We first note that 
\begin{equation}\label{pistar}\pi^*_T(K_T(\X))\subset K_T(\Y)^{W_G}.\end{equation}
This is because for any $w\in W_G=N_G(T)/T$ we have
$\pi\circ w =\pi$ (here $w$ is regarded as an automorphism of $\Y=Sp(n)/T$,
see equation (\ref{weyl})).

As mentioned in the introduction, the space
$\Y=Sp(n)/T$ is a  complete complex flag variety.
Thus, the fixed point set of the $T$ action is  
$$(Sp(n)/T)^T=W_{Sp(n)}.$$
The image of this set under $\pi$ is $Fl_n(\H)^T=\S_n$.
In fact, the restriction of $\pi$ to 
the fixed point set is the canonical projection
$W_{Sp(n)}\to W_{Sp(n)}/W_G =\S_n$
(see Section \ref{roots}).
The homomorphism between the $T$-equivariant 
$K$-theories induced by this map is the obvious map
$$p:\prod_{\S_n}R[T] \to \prod_{W_{Sp(n)}}R[T],$$
which is injective.
This homomorphism is the bottom arrow in the following commutative diagram
(where $\imath_T^*$ is the restriction homomorphism).

$$\xymatrix{
     K_T(\X) \ar[rr]^{\pi^*_T}   \ar@{->}[dd]^{i^*_T} & &
    \ar@{->}[dd]^{\imath_T^*}  K_T(\Y)  \\
     \\
   \ar[rr]^{{p}}  \prod_{\S_n}R[T] & & 
   \prod_{W_{Sp(n)}}R[T] } $$

\noindent The maps $i^*_T$ and $\imath^*_T$ are also injective, by Proposition
\ref{inject}, respectively the fact that $\Y$ is a complex flag variety
(for such spaces, injectivity is proved for instance in \cite[Lemma 2.2]{Mc}).
We deduce that $\pi^*_T$ is injective.


We now prove that the image of $\pi^*_T$ is the whole
$K_T(\Y)^{W_G}$ and that 
 $\overline{\kappa}_T$ is an
isomorphism.
Since $\Y=Sp(n)/T$ we deduce that  
\begin{equation}\label{isomo}K_T(\Y) = R[T] \otimes_{R[Sp(n)]} R[T].
\end{equation}
More precisely, by
 \cite[The Main Theorem]{Mc}, we have an isomorphism
\begin{equation}\label{isomor} R[T] \otimes_{R[Sp(n)]} R[T] \to K_T(\Y)\end{equation}
 which is the $R[T]$-linear extension of the map
 $\sigma \mapsto V'_\sigma$, for any complex representation 
 $\sigma : T\times \bC^m \to \bC^m$.
 Here we define 
\begin{equation}\label{gkvp}V'_\sigma=Sp(n)\times \bC^m/((k,v)\sim (kg,\sigma(g^{-1})v), \ 
{\rm for \ all} \ k\in Sp(n), g\in T, v\in \bC^m),\end{equation}
which is a complex vector bundle over $\Y=Sp(n)/T$, with an obvious
$T$ action.  
We claim that 
the composed homomorphism
$$\pi_T^* \circ \overline{\kappa}_T:R[T]\otimes_{R[Sp(n)]}R[G]
\to  R[T] \otimes_{R[Sp(n)]} R[T]$$ is induced by the canonical embedding 
$R[G] \rightarrow R[T]$: 
this follows  from equations (\ref{gkv}) and (\ref{gkvp}) and the
fact that $\pi_T^*$ maps a vector bundle to its  pull-back via $\pi_T$.
Consequently, $\pi^*_T\circ \overline{\kappa}_T$ is an isomorphism
between its domain and $R[T] \otimes_{R[Sp(n)]} R[G]$.
This space is just $K_T(\Y)^{W_G}$ (as we will see below). Thus,
 the image of $\pi^*_T$ is 
 $K_T(\Y)^{W_G}$ (see also equation (\ref{pistar})).
We also deduce that 
 $\overline{\kappa}_T$ is an isomorphism.
 
  It only remains to show that, via the identification (\ref{isomo}), we have
 $$K_T(\Y)^{W_G}=R[T]\otimes_{R[Sp(n)]}R[G].$$ 
 Indeed, this is a consequence of the fact that the isomorphism given by
 equation (\ref{isomor}) is $W_{Sp(n)}$-equivariant with respect to the  action
 given by $$w.(\chi_1\otimes \chi_2) =\chi_1\otimes (w\chi_2),$$
for all $w\in W_{Sp(n)}$, $\chi_1,\chi_2\in R[T]$ (see \cite[Section 1.4]{Mc} for a proof);
we deduce that
$$K_T(\Y)^{W_G} = R[T] \otimes_{R[Sp(n)]} R[T]^{W_G}
=R[T] \otimes_{R[Sp(n)]} R[G],$$
as required.
 The proposition is now proved.
\end{proof}

We now focus on the last statement in Theorem \ref{firstmain}.
Since $\Y$ is a complex flag variety, $K_T(\Y)$ can be identified with the Grothendieck group of the $T$-equivariant coherent sheaves on $\Y$, and the Bruhat decomposition gives an $R[T]$-basis of $K_T(\Y)$. More precisely, this  is the set of all classes $[\mathcal{O}_w]$ of the structure sheaves of the Schubert varieties $\Y_w$, for $w \in W_{Sp(n)}$
(cf. e.g. \cite[Section 4]{Ko-Ku}).

The following lemma will be needed later.

\begin{lemma}\label{mus}

Let $w \in W_{Sp(n)}$ and let $s_\nu$ be a simple reflection of $W_{Sp(n)}$.
 If $ws_{\nu} \le w$ in the Bruhat order, then $s_{\nu}[\mathcal{O}_w] = [\mathcal{O}_w]$.

\end{lemma}
\begin{proof}
The multiplication by $s_{\nu}$ from the right (see equation
(\ref{weyl})) induces a smooth map from $\Y$ to $\Y$. 
Consequently, since $\Y$ is a smooth variety, for any (not necessarily smooth) subvariety $Z$ of 
$Sp(n)/T$, we have $s_{\nu} [\mathcal{O}_Z] = 
[\mathcal{O}_{s_{\nu}(Z)}]$ 
(note that $s_{\nu}^{-1} = s_{\nu}$). In particular, for $w \in W_{Sp(n)}$ we have  
$s_{\nu} [\mathcal{O}_{w}] = [\mathcal{O}_{s_{\nu}(\Y_w)}]$. 
Since $ws_{\nu} \le  w$, we have  $s_{\nu}(\Y_w)=\Y_w$ (as we will see  below),
hence $s_{\nu} [\mathcal{O}_{w}]=[\mathcal{O}_{w}]$, as needed.

We still need to prove that $s_\nu(\Y_w)=\Y_w$. Since $s_\nu^{-1} = s_\nu$, it is enough to show that $s_\nu(\Y_w) \subset \Y_w$. 
We have $$\Y=Sp(n)/T=Sp(2n,\bC)/\B,$$
where $\B$ is a Borel subgroup of $Sp(2n,\bC)$ with $T\subset \B$. 
By definition, $\Y_w=\overline{\B w\B /\B}$.
Hence we have $$s_\nu(\Y_w)= s_\nu(\overline{\B w\B /\B}) \subset 
\overline{s_\nu(\B w\B /\B)}\subset
\overline{(\B w\B \B s_\nu \B)/\B}.$$ 
From the general theory of Tits systems  (see for example \cite{Bo}), we know that 
$\B w\B\B s_\nu \B \subset  \B w\B \cup \B ws_\nu \B$. 
This union is contained in $\Y_w$, since $ws_\nu\le  w$ in the Bruhat order.
The conclusion follows.
\end{proof}

From here we deduce the desired result (see again Theorem \ref{firstmain}),
namely:

\begin{proposition}\label{maxl}
The set $$\{ [\mathcal{O}_{w}] \ : \  w \in W_{Sp(n)} {\it  \ is \ a  \ maximal \ length \ representative \ in \, } W_{Sp(n)}/W_G \}$$ is a $R[T]$-basis of $\pi^*_T(K_T(\X))$.
\end{proposition}

\begin{proof} 
Let us denote by $W$ the subset of $W_{Sp(n)}$ consisting of all maximal length representatives of $W_{Sp(n)}/W_G$. 
The previous lemma implies that if $v\in W_G$ and 
$w\in W$, then $$v[\O_w]=[\O_w].$$ 
More precisely, we write $v=s_{1}\ldots s_{k}$, where
$s_{1},\ldots, s_{k}\in \{s_{2L^1},\ldots,s_{2L^n}\}$, which is the generating set 
of $W_G$ (see Section \ref{roots}); by Lemma \ref{mus} we have $s_\nu [\O_w]=[\O_w],$ for all $\nu\in \{1,2, \ldots, k\}$.  
We deduce from Proposition \ref{thretwo}  that  $\pi^*_T(K_T(\X))$ contains all
 $[\O_w]$ with $w\in W$.
 
To prove the converse inclusion, let us suppose that there exists  $ \xi \in K_T(\X)$ such that $\xi \notin \bigoplus_{w \in W} R[T][\mathcal{O}_w]$.  Since $K_T(\Y) = \bigoplus_{w \in W_{Sp(n)}} R[T][\mathcal{O}_w]$, we can write
$$\xi= \sum_{w \in W_{Sp(n)}} a_w [\mathcal{O}_w],$$
where $a_w \in R[T].$
We deduce $$\xi': = \xi -  \sum_{w \in W} a_w [\mathcal{O}_w] \in K_T(\X) \setminus  \{  0\}.$$ Consequently, the set $\{ \xi'\}\cup \{[\mathcal{O}_w] \ : \  w \in W \}$ is an
$R[T]$-free family of  $n! + 1$ elements in $K_T(\X)$. This is not possible,
since, by Proposition \ref{inject}, $K_T(\X)$ is a free $R[T]$-module of rank $n!$. The contradiction finishes the proof. 
\end{proof}

\noindent{\bf Remark.} We would like to point out that a result similar to
Proposition \ref{maxl} holds true for the usual (that is, non-equivariant)
$K$-theory group of $Fl_n(\H)$. 
Namely, from Pittie's theorem \cite{Pi}, we have
$K(\X)=R[G]/R[Sp(n)]$ and 
$K(\Y)=R[T]/R[Sp(n)]$. In terms of these identifications, the
homomorphism
 $\pi^*: K(\X)\to K(\Y)$ induced by $\pi:\Y\to \X$ is the inclusion induced by
the (injective) map $R[G]\to R[T]$ which assigns to a representation of
$G$ its restriction to $T$.
We deduce that $$K(\X)=K(\Y)^{W_G}.$$
Consequently, $K(\X)$ is the subring of $K(\Y)$ generated by
the elements of the Schubert basis induced by  $w\in W_{Sp(n)}$
which are maximal length representatives of $W_{Sp(n)}/W_G$:
this can be proved  by using the non-equivariant analogue of Proposition
\ref{maxl}.

We conclude the section with a GKM description of
$K_T(\X)$. We will deduce it from the GKM description
of $K_T(\Y)$ (recall that $\Y$ is a complete complex flag variety)
and the fact that $K_T(\X)=K_T(\Y)^{W_G}$ (see Proposition 
\ref{thretwo}). The notations established in Section \ref{roots} are
used in what follows.  

First, we  identify $\Y^T=W_{Sp(n)}$.
By  \cite[Theorem 1.6]{Mc},  the ring homomorphism 
$K_T(\Y)\to K_T(W_{Sp(n)})$
induced by the inclusion map is injective and its image is
$$   \{  (f_{{w}}) \in   \!\!\!\!\!  \!\!\! \prod_{{w} \in  W_{Sp(n)}/W_G} \!\!\!\!\!  \!\!\! R[T]  \ : \    e^{\alpha}-1  {\rm \,\,  divides \,\, }  f_{{w}}-f_{{s_{\alpha}w}} \,\,\, 
{\rm  for \ all \ }  \alpha \in \Delta^+   \}.$$
The isomorphism is $W_{Sp(n)}$-equivariant if we let the Weyl group $W_{Sp(n)}$ act on the space in the right-hand side of the previous equation
by 
$$v.(f_w) = (f_{wv^{-1}}) $$
for all $v\in W_{Sp(n)}$.
This can be proved as follows: take  $v,w\in W_{Sp(n)}$ and consider the
ring automorphism $v^*$ 
of $K_T(\Y)$ induced by $v$ (see equation (\ref{weyl})), as well as
the map $i_w^*:K_T(\Y)\to K_T(\{w\})$ induced by the inclusion
$\{w\}\hookrightarrow \Y$; for any $x\in K_T(\Y)$ we have
$$i_w^*(v^*(x)) =(v\circ i_w)^*(x)=i_{v(w)}^*(x)
=i_{wv^{-1}}^*(x).$$
Consequently,  $K_T(\Y)^{W_{G}}$ can be identified with 
$$   \{  (f_{\overline{w}}) \in   \!\!\!\!\!  \!\!\! \prod_{\overline{w} \in  W_{Sp(n)}/W_G} \!\!\!\!\!  \!\!\! R[T]  \ : \    e^{\alpha}-1  {\rm \,\,  divides \,\, }  f_{\overline{w}}-f_{\overline{s_{\alpha}w}} \,\,\, 
{\rm  for \ all \ }  \alpha \in \Delta^+   {\rm \,\, such \,\, that  \,\,} s_{\alpha} \notin W_{G} \}.$$
Here $ W_{Sp(n)}/W_G$ denotes the set of right cosets $\overline{w}=wW_G$ with 
$w\in W_{Sp(n)}$. The reason why in the previous equation we only need to consider 
roots $\alpha \in   \Delta^+$ such that  $s_{\alpha} \notin W_{G}$ is that
if $s_\alpha$ does belong to $W_G$ then 
$$\overline{s_{\alpha}w} = \overline{w w^{-1}s_{\alpha}w}=\overline{w},$$
for any $w\in W_{Sp(n)}$ (because $W_G$ is a normal subgroup of $W_{Sp(n)}$). 
The roots $\alpha \in \Delta^+$ such that  $s_\alpha \notin W_G$ are 
$L^\mu-L^\nu$ and $L^\mu+L^\nu$, where $1\le \mu <\nu \le n$. 
We saw in  Section \ref{roots} that via  the identification
$ W_{Sp(n)}/W_G=\S_n$,
we have $s_{L^\mu-L^\nu}W_G=s_{L^\mu+L^\nu}W_G=(\mu,\nu)$.
Let us denote $$e^{L_{\nu}}=x_{\nu}, \ 1\le \nu \le n.$$
This implies that
$e^{L^\mu-L^\nu}-1=x_\mu x_\nu^{-1}-1$ and $e^{L^\mu+L^\nu}-1=x_\mu 
x_\nu-1.$
Since these two polynomials are relatively prime in $\mathbb{Z}[x_1^{\pm 1}, \ldots, 
x_n^{\pm 1} ]$, we deduce the following proposition. 
 \begin{proposition}\label{gkmc} The  homomorphism
 $K_T(Fl_n(\H))\to K_T(Fl_n(\H)^T)$ induced by the inclusion
 $Fl_n(\H)^T\to Fl_n(\H)$ is injective.  
   Its image  is equal to
      $$\{  (f_\tau) \in    \!\!\! \prod_{\tau \in \S_n} \!\!\! \mathbb{Z}[x_1^{\pm 1}, \ldots, 
x_n^{\pm 1} ]  \ : \   
(x_\mu x_\nu^{-1}-1)(x_\mu x_\nu -1) \ {\it divides } \ f_\tau-f_{(\mu,\nu)\tau}  \,\,\, {\it for \ all} \   1 \leq \mu < \nu \leq n \}.$$
\end{proposition}

\section{$G$-equivariant $K$-theory}\label{lasts}

In this section we will prove Theorem \ref{secondmain}.
Like in the previous section, we denote $\X=Fl_n(\H)$,
which is the same as the homogeneous space $Sp(n)/G$.
We consider the following commutative diagram:
$$\xymatrix{
     K_G(\X) \ar[rr]^{j^*}   \ar@{->}[dd]^{i_{G}^*} & &
    \ar@{->}[dd]^{i_{T}^*}  K_T(\X)  \\
     \\
   \ar[rr]^{\tilde{p}}  \prod_{\S_n}R[G] & & 
   \prod_{\S_n}R[T] } $$
Here $j^*$ and $\tilde{p}$ are induced by the restriction of
the $G$ action to $T$. The maps 
$i_T^*$, $\tilde{p}$, and $j^*$ are injective:
the first by Proposition \ref{gkmc}, the second by
e.g. \cite[Chapter 13, Section 8]{Hu}, and the last  by \cite[Theorem 4.4]{Mc}. 
We deduce that $i^*_G$ is injective too.

 We are interested in the image of $i^*_G$.
Let us consider  
 the action of $W_G$ on $R[T] \otimes_{R[Sp(n)]} R[G]$
 given by
 $$w.(\chi_1\otimes \chi_2)=(w\chi_1) \otimes \chi_2,$$
for any $\chi_1\in R[T], \chi_2\in R[G], w\in W_{G}$,
as well as the 
diagonal action of $W_G$ on
 $\prod_{\S_n}R[T]$. We need the following lemma.
 \begin{lemma}\label{themap} The map
 $$i^*_T\circ \overline{\kappa}_T : 
 R[T] \otimes_{R[Sp(n)]} R[G]\to  \prod_{\S_n}R[T]$$
 is a $W_G$-equivariant homomorphism.
 \end{lemma}
 \begin{proof} We take $w\in \X^T$ 
 and  show that the map
 $ R[T] \times R[G]\to  R[T]$ given by
\begin{equation}\label{chi} (\chi_1,\chi_2)\mapsto i^*_T(
\kappa_T(\chi_1\otimes\chi_2))_w
\end{equation}
 is $W_G$-equivariant (the map $\kappa_T$ was defined 
 in Section \ref{prev}). We identify $\X=Sp(n)/G$ and 
 write $w=kG$, where $k\in N_{Sp(n)}(G)$. 
 Assume that $\chi_1,\chi_2$ in  equation (\ref{chi}) are the characters
 of the representations $(V_1,\rho_1)$, respectively
 $(V_2,\rho_2)$. An easy exercise  shows that 
 $\kappa_T(\chi_1\otimes \chi_2)_w$ is the $T$-representation
 on $V_1\otimes V_2$ given by
 $$t.(v_1\otimes v_2) = v_1 \otimes \rho_2(k^{-1}tk)(v_2),$$
 for all $v_1\in V_1, v_2\in V_2$ and $t\in T$. This is the tensor
 product of two $T$-representations, the first one being trivial and
 the second one lying in
 $R[G]$, which is the same as $R[T]^{W_G}$.
 The $W_G$-equivariance of our map is now clear.  
   \end{proof}
 
From now on we will identify
$$K_T(\X)=R[T] \otimes_{R[Sp(n)]} R[G],$$
by using Proposition \ref{thretwo}. In this way, $W_G$ 
acts on $K_T(\X)$. From the previous lemma we
deduce that the map 
 $i^*_T$ is a $W_G$-equivariant homomorphism.
  Since $R[G] = R[T]^{W_{G}}$, the image of the map
  $$\tilde{p} \circ  i_{G}^*=i_T^*\circ j^*$$ is included in 
 $(i_T^*(K_T(\X)))^{W_{G}}$.
 Since $i^*_T$ is injective, we deduce that
 the image of $j^*$ is contained in $K_T(\X)^{W_G}$.

Like in the $T$-equivariant case (see the previous  section), we consider the homomorphism 
$$\overline{\kappa}_G : R[G] \otimes_{R[Sp(n)]} R[G] \rightarrow K_G(\X),$$
which is $R[G]$-linear and satisfies
$$\overline{\kappa}_G(1\otimes \chi)=[V_\rho]$$
for any representation $\rho: G\times V \to V$ of character $\chi$.
For the definition of $V_\rho$, see equation (\ref{vro}).  
The composition $$j^* \circ \overline{\kappa}_G  : R[G] \otimes_{R[Sp(n)]} R[G] \rightarrow K_T(\X) =  R[T] \otimes_{R[Sp(n)]} R[G]$$ is induced by the restriction map from $G$ to $T$. Since $R[G] = R[T]^{W_{G}}$, the map $j^* \circ \overline{\kappa}_G$ is an isomorphism between $R[G] \otimes_{R[Sp(n)]} R[G]$ and its image in $K_T(\X)$.
The image is actually equal to
$K_T(\X)^{W_{G}}$.
Consequently, the image of $j^*$ is $K_T(\X)^{W_{G}}$. 
Since $j^*$ is injective (as we saw above), we deduce that   it is an
isomorphism between $K_G(\X)$ and $K_T(\X)^{W_G}$.
Incidentally, we have also proved the following result:
\begin{proposition}\label{themap1} The map $\overline{\kappa}_G$ defined above is a 
ring isomorphism.
\end{proposition}

From the commutative diagram at the beginning of the section we deduce  that the image of $i^*_G$  
consists of all $W_G$-invariants in the image of $i^*_T$.
By Proposition \ref{gkmc}, the image of $i^*_T$ is
$$\{  (f_\tau) \in    \!\!\! \prod_{\tau \in \S_n} \!\!\! \mathbb{Z}[x_1^{\pm 1}, \ldots, 
x_n^{\pm 1} ] ^{W_{G}} \ : \   
(x_\mu x_\nu^{-1}-1)(x_\mu x_\nu -1) \ {\rm divides } \ f_\tau-f_{(\mu,\nu)\tau}  \,\,\, {\rm for \ all} \   1 \leq \mu < \nu \leq n \}.$$
By  Lemma \ref{last}, below, this is the same as
$$\{  (f_\tau) \in    \!\!\! \prod_{\tau \in \S_n} \!\!\! \mathbb{Z}[X_1, \ldots, X_n ]  \ : \  
f_\tau-f_{ {(\mu,\nu)}\tau}  {\rm \,\, is \,\, divisible \,\, by \,\,} X_\mu-X_\nu  \,\,\, {\rm for \ all} \   1 \leq \mu < \nu \leq n \}.
$$
Theorem \ref{secondmain} is now completely proved.
The following  lemma has been used above.

\begin{lemma}\label{last} Let $x_1,\ldots, x_n$ be some variables and set
$$X_\nu=x_\nu+x_\nu^{-1},$$
$1 \le \nu \le n$. 
An element of $\bZ[X_1,\ldots, X_n]$ is divisible by
$(x_\mu x_\nu^{-1} -1)(x_\mu x_\nu -1)$ in the ring
$\bZ[x_1^{\pm 1},\ldots,x_n^{\pm 1}]$ 
if and only if it is divisible by $X_\mu -X_\nu$
in $\bZ[X_1,\ldots,X_n]$.
\end{lemma}

\begin{proof} We first prove the following claim.

\noindent {\it Claim.} The ring $\bZ[x_1^{\pm 1},\ldots,x_n^{\pm 1}]$ 
is a free module over $\bZ[X_1,\ldots, X_n]$ of basis
$x_1^{\epsilon_1}\ldots x_n^{\epsilon_n}$, where
$\epsilon_\nu\in \{-1,0\}$, $1\le \nu \le n$.

We first notice that if $x$ is a variable,
${ X}=x+x^{-1}$, and  $R$ is an arbitrary unit ring, then
$R[x, x^{-1}]$ is a free module over $R[{ X}]$
of basis $1, x^{-1}$. 
This implies the claim by a recursive argument.
Namely, we take successively 
 $R=\bZ[x_{\mu}^{\pm 1},\ldots,x_n^{\pm 1}]$
 for $\mu=2, \ldots, n$, which gives $R[x_{\mu-1}^{\pm 1}]=\bZ[x_{\mu-1}^{\pm 1},\ldots,x_n^{\pm 1}]$.

We now turn to the proof of the lemma. For sake of simplicity let us make
$\mu=1$ and $\nu=2$. We   note that
$$X_1-X_2= x_1^{-1} (x_1x_2^{-1}-1)(x_1x_2-1).$$
Thus, if $f$ is divisible by $X_1-X_2$ 
then it is divisible by
$(x_1 x_2^{-1} -1)(x_1x_2 -1)$.
We will now prove the converse.
Assume that  $f\in \bZ[X_1,\ldots, X_n]$ 
is divisible by $(x_1x_2^{-1}-1)(x_1x_2-1)$.
We deduce that 
\begin{equation}\label{laste}f=(X_1-X_2)g,\end{equation}
where $g\in \bZ[ x_1^{\pm 1},\ldots,x_n^{\pm 1}]$.
We consider the expansion of $g$ with respect to the basis indicated in the claim
and denote by $g_0\in \bZ[X_1,\ldots, X_n]$ the coefficient of 1. Equation (\ref{laste})
implies $$f=(X_1-X_2)g_0.$$
This finishes the proof.
\end{proof}

\bigskip
\bigskip


 \bibliographystyle{plain}

\end{document}